\documentclass[11pt]{article}

\usepackage{amscd,amstex,amssymb}

\def\eqref#1{(\ref{#1})}
\newcommand{\goth}{\mathfrak}
\newcommand{\arrow}{{\:\longrightarrow\:}}

\newcommand{\C}{{\Bbb C}}
\newcommand{\R}{{\Bbb R}}

\newcommand{\inangles}[1]{{\langle #1\rangle}}

\newcommand{\restrict}[1]{{\left|_{{\phantom{|}\!\!}_{#1}}\right.}}

\renewcommand{\c}[1]{{\cal #1}}
\newcommand{\calo}{{\cal O}}

\renewcommand{\tilde}{\widetilde}

\renewcommand{\phi}{\varphi}
\renewcommand{\epsilon}{\varepsilon}

\newcommand{\End}{\operatorname{End}}

\newcommand{\Def}{\operatorname{Def}}

\newcommand{\comment}[1]{{}}

\def\blacksquare{\hbox{\vrule width 4pt height 4pt depth 0pt}}
\def\endproof{\blacksquare}

\makeatletter

\@ifundefined{Bbb}
     {\newcommand{\Bbb}[1]{{\mathbb #1}}}%
{}%     {\edef\Bbb#1{{\Bbb #1}}}

\newcommand{\ps@verbit}{%
  \renewcommand{\@oddhead}{%
          \scriptsize
          {Desingularization of hyperk\"ahler varieties II}
          \hfil\tiny {final version, December 12, 1996}}
  \renewcommand{\@evenhead}{\@oddhead}
  \renewcommand{\@oddfoot}{\hfil\thepage\hfil}
  \renewcommand{\@evenfoot}{\@oddfoot}}
 
\newcounter{Mycounter}[section]
\newcounter{lemma}[section]
\renewcommand{\thelemma}{{Lemma \thesection.\arabic{lemma}}}
\newcommand{\lemma}{%
     \setcounter{lemma}{\value{Mycounter}}
     \refstepcounter{lemma}
     \stepcounter{Mycounter}
     {\bf \thelemma:\ }}

\newcounter{claim}[section]
\renewcommand{\theclaim}{{Claim \thesection.\arabic{claim}}}
\newcommand{\claim}{%
     \setcounter{claim}{\value{Mycounter}}
     \refstepcounter{claim}
     \stepcounter{Mycounter}
     {\bf \theclaim:\ }}

\newcounter{sublemma}[section]
\newcounter{corollary}[section]
\renewcommand{\thecorollary}{{Corollary \thesection.\arabic{corollary}}}
\newcommand{\corollary}{%
     \setcounter{corollary}{\value{Mycounter}}
     \refstepcounter{corollary}
     \stepcounter{Mycounter}
     {\bf \thecorollary:\ }}

\newcounter{theorem}[section]
\renewcommand{\thetheorem}{{Theorem \thesection.\arabic{theorem}}}
\newcommand{\theorem}{%
     \setcounter{theorem}{\value{Mycounter}}
     \refstepcounter{theorem}
     \stepcounter{Mycounter}
     {\bf \thetheorem:\ }}

\newcounter{conjecture}[section]
\newcounter{proposition}[section]
\renewcommand{\theproposition}
       {{Proposition \thesection.\arabic{proposition}}}
\newcommand{\proposition}{%
     \setcounter{proposition}{\value{Mycounter}}
     \refstepcounter{proposition}
     \stepcounter{Mycounter}
     {\bf \theproposition:\ }}

\newcounter{definition}[section]
\renewcommand{\thedefinition}
       {{Definition \thesection.\arabic{definition}}}
\newcommand{\definition}{%
     \setcounter{definition}{\value{Mycounter}}
     \refstepcounter{definition}
     \stepcounter{Mycounter}
     {\bf \thedefinition:\ }}

\newcounter{example}[section]
\newcounter{remark}[section]
\renewcommand{\theremark}{{Remark \thesection.\arabic{remark}}}
\newcommand{\remark}{%
     \setcounter{remark}{\value{Mycounter}}
     \refstepcounter{remark}
     \stepcounter{Mycounter}
     {\bf \theremark:\ }}

\newcounter{problem}[section]
\newcounter{question}[section]
\@addtoreset{equation}{section}
\@addtoreset{footnote}{section}
\makeatother

\begin{document}

\begin{center}
{\Large\bf
	Desingularization of singular hyperk\"ahler varieties II.}\\[4mm]
Misha Verbitsky,\footnote{Supported by the NSF grant 9304580}\\[4mm]
{\tt verbit@@thelema.dnttm.rssi.ru, verbit@@math.ias.edu}
\end{center}

\hfill

{\small 
\hspace{0.2\linewidth}
\begin{minipage}[t]{0.7\linewidth}
{\bf Abstract:} We construct 
a natural hyperk\"ahler
desingularization for all singular hyperk\"ahler 
varieties. 
\end{minipage}
}

\tableofcontents

%%%%%%%%%%%%%%%%%%%%%%%%%%%%%%%%%%%%%%%%

\section{Introduction}

%%%%%%%%%%%%%%%%%%%%%%%%%%%%%%%%%%%%%%%%

A hyperk\"ahler manifold is a Riemannian manifold with an action
of a quaternion algebra $\Bbb H$ in its tangent bundle, such that
for all $I\in \Bbb H$, $I^2=-1$, $I$ establishes a complex,
K\"ahler structure on $M$
(see \ref{_hyperkahler_manifold_Definition_} for details).
We extend this definition to singular varieties.
The notion of a singular hyperk\"ahler variety has its
origin in \cite{_Verbitsky:Hyperholo_bundles_} 
(see also \cite{_Verbitsky:Deforma_} and 
\cite{_Verbitsky:Desingu_}). Examples of singular
hyperk\"ahler varieties are numerous, and come from several
diverse sources (\ref{_singu_hype_Remark_},
\ref{_hyperho_defo_hyperka_Theorem_};
for additional examples see \cite{_Verbitsky:Deforma_},
Section 10). There is a weaker version of this definition:
a notion of {\bf hypercomplex variety}. Singular hypercomplex
varieties is what this paper primarily deals with.

For a real analytic variety $M$, we say that $M$ is {\bf hypercomplex},
if $M$ is equipped with the complex structures $I$, $J$ and $K$,
such that $I\circ J = - J\circ I =K$, and certain integrability
conditions are satisfied (for a precise statement, see 
\ref{_hypercomplex_Definition_}). The paper 
\cite{_Verbitsky:Desingu_} dealt with the hypercomplex varieties
with ``locally homogeneous singularities'' (LHS). A complex analytic
variety $M$ is LHS if for each point $x\in M$, the completion
$A$ of the local ring $\calo_x M$ can be represented as 
a quotient of the power series ring by a homogeneous ideal
(\ref{_SLHS_Definition_}, \ref{_locally_homo_coord_Claim_}).
In \cite{_Verbitsky:Desingu_}, we described explicitly
the singularities of hypercomplex LHS varieties.
We have shown that every such variety, considered as a complex
variety with a complex structure induced 
from the quaternions, is locally isomorphic
to a union of planes in $\C^n$ (\ref{_singula_stru_Theorem_}).
The normalization of such a variety is non-singular,
which follows from this description of singularities.
This gives a canonical, functorial way to desingularize hyperk\"ahler
and hypercomplex varieties (\ref{_desingu_Theorem_}).

The purpose of the present paper is to show that all
hypercomplex varieties have locally homogeneous singularities 
(\ref{_hyperco_SLHS_Theorem_}).
This is used to extend the desingularization results to all
singular hyperk\"ahler or hypercomplex varieties 
(\ref{_hyperco_desingu_Corollary_}).

In Sections 
\ref{_comple_with_au_Section_}--\ref{_homogeni_on_hype_Section_},
we prove that all hypercomplex varieties have locally homogeneous 
singularities. Section \ref{_comple_with_au_Section_} is purely
a commutative algebra. We work with a complete local Noetherian
ring $A$ over $\C$. By definition, an automorphism $e:\; A \arrow A$
is called {\bf homogenizing} (\ref{_homogeni_automo_Definition_})
if its differential acts as
a dilatation on the Zariski tangent space of $A$, with dilatation
coefficient $|\lambda|<1$. As usual, by the Zariski tangent space
we understand the space $\goth m_A /\goth m_A^2$, where $\goth m_A$
is a maximal ideal of $A$.

The main result of Section \ref{_comple_with_au_Section_} is the
following. For a complete local Noetherian
ring $A$ over $\C$ equipped with a homogenizing automorphism
$e:\; A \arrow A$, we show that $A$ has locally homogeneous 
singularities.

In Section \ref{_homogeni_on_hype_Section_}, we construct a 
natural homogenizing automorphism of the ring of germs of
complex analytic functions on a hypercomplex variety $M$
(\ref{_homogenizing_Proposition_}). Applying Section 
\ref{_comple_with_au_Section_}, we obtain that every 
hypercomplex variety has locally homogeneous 
singularities.

%%%%%%%%%%%%%%%%%%%%%%%%%%%%%%%%%%%%%%%%

\section{Preliminaries}

%%%%%%%%%%%%%%%%%%%%%%%%%%%%%%%%%%%%%%%%

\subsection{Definitions}
 
This subsection contains a compression of 
the basic definitions from hyperk\"ahler geometry, found, for instance, in
\cite{_Besse:Einst_Manifo_} or in \cite{_Beauville_}.
 
\hfill
 
%%%%%%%%%%%%%%%%%%%%%%%%%%%%%%%%%%%%%%%%%%%%%%%%%%%%%%%%%%%%%%%%%
\definition \label{_hyperkahler_manifold_Definition_} %%%%%%%%%%
(\cite{_Besse:Einst_Manifo_}) A {\bf hyperk\"ahler manifold} is a
Riemannian manifold $M$ endowed with three complex structures $I$, $J$
and $K$, such that the following holds.
 
\begin{description}
\item[(i)]  the metric on $M$ is K\"ahler with respect to these complex 
structures and
 
\item[(ii)] $I$, $J$ and $K$, considered as  endomorphisms
of a real tangent bundle, satisfy the relation 
$I\circ J=-J\circ I = K$.
\end{description}
 
\hfill 
 
The notion of a hyperk\"ahler manifold was 
introduced by E. Calabi (\cite{_Calabi_}).

\hfill
 
Clearly, a hyperk\"ahler manifold has the natural action of
quaternion algebra ${\Bbb H}$ in its real tangent bundle $TM$. 
Therefore its complex dimension is even.
For each quaternion $L\in \Bbb H$, $L^2=-1$,
the corresponding automorphism of $TM$ is an almost complex
structure. It is easy to check that this almost 
complex structure is integrable (\cite{_Besse:Einst_Manifo_}).
 
\hfill
 
%%%%%%%%%%%%%%%%%%%%%%%%%%%%%%%%%%%%%%%%%%%%%%%%%%%%%%%%%%%%
\definition \label{_indu_comple_str_Definition_} 
Let $M$ be a hyperk\"ahler manifold, $L$ a quaternion satisfying
$L^2=-1$. The corresponding complex structure on $M$ is called
{\bf an induced complex structure}. The $M$ considered as a complex
manifold is denoted by $(M, L)$.
 
\hfill
 
Let $M$ be a hyperk\"ahler manifold. We identify the group $SU(2)$
with the group of unitary quaternions. This gives a canonical 
action of $SU(2)$ on the tangent bundle and all its tensor
powers. In particular, we obtain a natural action of $SU(2)$
on the bundle of differential forms. 

\hfill

%%%%%%%%%%%%%%%%%%%%%%%%%%%%%%%%%%%%%%%%%%%%%%%%%%%%%%%%%%%%
\lemma \label{_SU(2)_commu_Laplace_Lemma_}
The action of $SU(2)$ on differential forms commutes
with the Laplacian.
 
{\bf Proof:} This is Proposition 1.1
of \cite{_Verbitsky:Hyperholo_bundles_}. \endproof
 
Thus, for compact $M$, we may speak of the natural action of
$SU(2)$ in cohomology.

%%%%%%%%%%%%%%%%%%%%%%%%%%%%%%%%%%%%%%%%%%%%%%%%
 
\subsection{Trianalytic subvarieties in compact hyperk\"ahler
manifolds.}
 
%%%%%%%%%%%%%%%%%%%%%%%%%%%%%%%%%%%%%%%%%%%%%%%%
 
In this subsection, we give a definition and a few basic properties
of trianalytic subvarieties of hyperk\"ahler manifolds. 
We follow \cite{_Verbitsky:Symplectic_II_}.
 
\hfill
 
Let $M$ be a compact hyperk\"ahler manifold, $\dim_\R M =2m$.
 
\hfill
 
%%%%%%%%%%%%%%%%%%%%%%%%%%%%%%%%%%%%%%%%%%%%%%%%%
\definition\label{_trianalytic_Definition_} %%%%%%%%%%
Let $N\subset M$ be a closed subset of $M$. Then $N$ is
called {\bf trianalytic} if $N$ is a complex analytic subset 
of $(M,L)$ for any induced complex structure $L$.
 
\hfill
 
Let $I$ be an induced complex structure on $M$,
and $N\subset(M,I)$ be
a closed analytic subvariety of $(M,I)$, $dim_\C N= n$.
Denote by $[N]\in H_{2n}(M)$ the homology class 
represented by $N$. Let $\inangles N\in H^{2m-2n}(M)$ denote 
the Poincare dual cohomology class. Recall that
the hyperk\"ahler structure induces the action of 
the group $SU(2)$ on the space $H^{2m-2n}(M)$.
 
\hfill
 
%%%%%%%%%%%%%%%%%%%%%%%%%%%%%%%%%%%%%%%%%%%%%%%%%%%%%%%%%%%%%%%%%
\theorem\label{_G_M_invariant_implies_trianalytic_Theorem_} %%%%%
Assume that $\inangles N\in  H^{2m-2n}(M)$ is invariant with respect
to the action of $SU(2)$ on $H^{2m-2n}(M)$. Then $N$ is trianalytic.
 
{\bf Proof:} This is Theorem 4.1 of 
\cite{_Verbitsky:Symplectic_II_}.
\endproof
 
\hfill

%%%%%%%%%%%%%%%%%%%%%%%%%%%%%%%%%%%%%%%%%%%%%%%%
\remark \label{_triana_dim_div_4_Remark_}
Trianalytic subvarieties have an action of quaternion algebra in
the tangent bundle. In particular,
the real dimension of such subvarieties is divisible by 4.
The non-singular part of a trianalytic subvariety is hyperk\"ahler.

%%%%%%%%%%%%%%%%%%%%%%%%%%%%%%%%%%%%%%%%%%%%%%%%%%
 
\subsection{Hypercomplex varieties}

%%%%%%%%%%%%%%%%%%%%%%%%%%%%%%%%%%%%%%%%%%%%%%%%%%

This subsection is based on the results and definitions from
\cite{_Verbitsky:Desingu_}.

Let $X$ be a complex variety, $X_\R$ the underlying real analytic variety.
In \cite{_Verbitsky:Desingu_}, Section 2, we constructed a natural 
automorphism of the sheaf of K\"ahler differentials on $X_\R$ 
\[ I:\; \Omega^1X_\R \arrow \Omega^1X_\R, \ \ I^2=-1. \]
This endomorphism is a generalization of the usual
notion of a complex structure operator, and its construction
is straightforward. We called $I$ {\bf the complex structure
operator on $X_\R$}. The operator $I$ is functorial:
for a morphism $f:\; X \arrow Y$ of complex varieties, 
the natural pullback map 
$df:\; f^*\Omega^1_{Y_\R} \arrow \Omega^1X_\R$
commutes with the complex structure operators 
(see \cite{_Verbitsky:Desingu_} for details).
The converse statement is also true:

\hfill

%%%%%%%%%%%%%%%%%%%%%%%%%%%%%%%%%%%%%%%%%%%%%%%%%%%%%%%%
\theorem \label{_commu_w_comple_str_Theorem_} %%%%%%%%%%
Let $X$, $Y$ be complex analytic varieties, and 
\[ f_\R:\; X_\R\arrow Y_\R\] be a morphism of underlying real
analytic varieties which commutes with the complex structure.
Then there exist a unique morphism $f:\; X\arrow Y$ of 
complex analytic varieties, such that $f_\R$ 
is its underlying morphism.

{\bf Proof:} This is \cite{_Verbitsky:Desingu_}, Theorem 2.1.
\endproof

\hfill

%%%%%%%%%%%%%%%%%%%%%%%%%%%%%%%%%%%%%%%%
\definition %%%%%%%%%%%%%%%%%%%%%%%%%%%%
Let $M$ be a real analytic variety, and
\[ I:\; \Omega^1(\calo_M)\arrow\Omega^1(\calo_M) \]
be an endomorphism satisfying $I^2=-1$. Then
$I$ is called {\bf an almost complex structure
on $M$}. If there exist a complex analytic structure $\goth C$
on $M$ such that $I$ appears as the complex structure operator
associated with $\goth C$, we say that $I$ is {\bf integrable}. 
\ref{_commu_w_comple_str_Theorem_} implies
that this complex structure is unique if it
exists. 

\hfill

%%%%%%%%%%%%%%%%%%%%%%%%%%%%%%%%%%%%%%%%%%%%%%%%%%
\definition \label{_hypercomplex_Definition_}
(Hypercomplex variety)
Let $M$ be a real analytic variety equipped with almost
complex structures $I$, $J$ and $K$, such that
$I\circ J = -J \circ I = K$. Assume that for all
$a, b, c\in \R$, such that $a^2 + b^2 + c^2=1$,
the almost complex structure $a I + b J + c K$ is integrable.
Then $M$ is called {\bf a hypercomplex variety}.

\hfill

%%%%%%%%%%%%%%%%%%%%%%%%%%%%%%%%%%%%%%%%
\remark
As follows from \cite{_Verbitsky:Desingu_}, Claim 2.7,
every hyperk\"ahler manifold is hypercomplex, in a natural way.
The proof is straightforward.
 
%%%%%%%%%%%%%%%%%%%%%%%%%%%%%%%%%%%%%%%%%%%%%%%%%%

\subsection{Singular hyperk\"ahler varieties}

%%%%%%%%%%%%%%%%%%%%%%%%%%%%%%%%%%%%%%%%%%%%%%%%%%

Throughout this paper, we never use the notion of
hyperk\"ahler variety. For our present purposes, the hypercomplex varieties
suffice. However, for the reader's benefit, we give a definition and
a list of examples of hyperk\"ahler varieties. All hyperk\"ahler
varieties are hypercomplex, and the converse is (most likely) false. 
However, it is difficult to construct examples of hypercomplex 
varieties which are not hyperk\"ahler, and all ``naturally'' occuring 
hypercomplex varieties come equipped with a singular hyperk\"ahler
structure.

This subsection is based on the results and definitions from
\cite{_Verbitsky:Deforma_} and \cite{_Verbitsky:Desingu_}.
For a more detailed exposition, the reader is referred to 
\cite{_Verbitsky:Deforma_}, Section 10.

\hfill

%%%%%%%%%%%%%%%%%%%%%%%%%%%%%%%%%%%%%%%%%%%%%%%%%%
\definition\label{_singu_hype_Definition_}
(\cite{_Verbitsky:Hyperholo_bundles_}, Definition 6.5)
Let $M$ be a hypercomplex variety (\ref{_hypercomplex_Definition_}).
The following data define a structure of a {\bf hyperk\"ahler variety}
on $M$.

\begin{description}

\item[(i)] For every $x\in M$, we have an $\R$-linear 
symmetric positively defined
bilinear form $s_x:\; T_x M \times T_x M \arrow \R$
on the corresponding real Zariski tangent space.

\item[(ii)] For each triple of induced complex structures
$I$, $J$, $K$, such that $I\circ J = K$, we have a 
holomorphic differential 2-form $\Omega\in \Omega^2(M, I)$.

\item[(iii)] 
Fix a triple of induced complex structure 
$I$, $J$, $K$, such that $I\circ J = K$. Consider the
corresponding differential 2-form $\Omega$ of (ii).
Let $J:\; T_x M \arrow T_x M$ be an endomorphism of 
the real Zariski tangent spaces defined by $J$, and $Re\Omega\restrict x$
the real part of $\Omega$, considered as a bilinear form on $T_x M$.
Let $r_x$ be a bilinear form $r_x:\; T_x M \times T_x M \arrow \R$ 
defined by $r_x(a,b) = - Re\Omega\restrict x (a, J(b))$.
Then $r_x$ is equal to the form $s_x$ of (i). In particular,
$r_x$ is independent from the choice of $I$, $J$, $K$.

\end{description}

%%%%%%%%%%%%%%%%%%%%%%%%%%%%%%%%%%%%%%%%%%%%%%%%%%
\noindent \remark \label{_singu_hype_Remark_}
\nopagebreak
\begin{description}
\item[(a)] It is clear how to define a morphism of hyperk\"ahler varieties.

\item[(b)]
For $M$ non-singular,  \ref{_singu_hype_Definition_} is
 equivalent to the usual
one (\ref{_hyperkahler_manifold_Definition_}). 
If $M$ is non-singular,
the form $s_x$ becomes the usual Riemann form, and 
$\Omega$ becomes the standard holomorphically symplectic form.

\item[(c)] It is easy to check the following.
Let $X$ be a hypercomplex subvariety of a hyperk\"ahler
variety $M$. Then, restricting the forms $s_x$ and $\Omega$
to $X$, we obtain a hyperk\"ahler structure on $X$. In particular,
trianalytic subvarieties of hyperk\"ahler manifolds are always
hyperk\"ahler, in the sense of \ref{_singu_hype_Definition_}.

\end{description}

\hfill

{\bf Caution:} Not everything which is seemingly hyperk\"ahler
satisfies the conditions of \ref{_singu_hype_Definition_}.
Take a quotient $M/G$ os a hyperk\"ahler manifold by an action 
of finite group $G$, acting in accordance with hyperk\"ahler
structure. Then $M/G$ is, generally speaking, {\it not} hyperk\"ahler
(see \cite{_Verbitsky:Deforma_}, Section 10 for details).

\hfill

The following theorem, proven in 
\cite{_Verbitsky:Hyperholo_bundles_} (Theorem 6.3), 
gives a convenient way to construct
examples of hyperk\"ahler varieties.

\hfill

%%%%%%%%%%%%%%%%%%%%%%%%%%%%%%%%%%%%%%%%%%%%%%%%%%
\theorem \label{_hyperho_defo_hyperka_Theorem_}
Let $M$ be a compact hyperk\"ahler manifold, $I$ an induced
complex structure and $B$ a stable holomorphic bundle over $(M, I)$.
Let $\Def(B)$ be the reduction\footnote{The deformation space might have
nilpotents in the structure sheaf. We take its reduction to avoid
this.} of the
deformation space of stable holomorphic structures on $B$.
Assume that $c_1(B)$, $c_2(B)$ are $SU(2)$-invariant, with respect
to the standard action of $SU(2)$ on $H^*(M)$. Then $\Def(B)$ has a
natural structure of a hyperk\"ahler variety. 

\nopagebreak
\endproof

%%%%%%%%%%%%%%%%%%%%%%%%%%%%%%%%%%%%%%%%

\section{Desingularization of hyperk\"ahler varieties}

%%%%%%%%%%%%%%%%%%%%%%%%%%%%%%%%%%%%%%%%

In this section, we recall the desingularization theorem
for the hypercomplex varieties with locally homogeneous singularities,
as it was proven 
in \cite{_Verbitsky:Desingu_}. In the last subsection,
we state the main result of this paper, which
is proven in Sections 
\ref{_comple_with_au_Section_}--\ref{_homogeni_on_hype_Section_}.

%%%%%%%%%%%%%%%%%%%%%%%%%%%%%%%%%%%%%%%%%%%%%%%%%%%%%%%%%%%%%%%%%%%%%%
\subsection{Spaces with locally homogeneous singularities.}

\noindent
%%%%%%%%%%%%%%%%%%%%%%%%%%%%%%%%%%%%%%%%%%%%%%%%%%
\definition
(local rings with LHS)
Let $A$ be a local ring. Denote by $\goth m$ its maximal ideal.
Let $A_{gr}$ be the corresponding associated graded ring.
Let $\hat A$, $\widehat{A_{gr}}$ be the $\goth m$-adic completion
of $A$, $A_{gr}$. Let $(\hat A)_{gr}$, $(\widehat{A_{gr}})_{gr}$ 
be the associated graded rings, which are naturally isomorphic to
$A_{gr}$. We say that $A$ {\bf has locally homogeneous singularities}
(LHS)
if there exists an isomorphism $\rho:\; \hat A \arrow \widehat{A_{gr}}$
which induces the standard isomorphism 
$i:\; (\hat A)_{gr}\arrow (\widehat{A_{gr}})_{gr}$ on associated
graded rings.

\hfill

%%%%%%%%%%%%%%%%%%%%%%%%%%%%%%%%%%%%%%%%%%%%%%%%%%
\definition\label{_SLHS_Definition_}
(SLHS)
Let $X$ be a complex or real analytic space. Then 
$X$ is called {be a space with locally homogeneous singularities}
(SLHS) if for each $x\in M$, the local ring $\calo_x M$ 
has locally homogeneous singularities.

\hfill
 
The following claim might shed a light on the origin of the term
``locally homogeneous singularities''.

\hfill

%%%%%%%%%%%%%%%%%%%%%%%%%%%%%%%%%%%%%%%%
\claim \label{_locally_homo_coord_Claim_}
Let $A$ be a complete local Noetherian ring over $\C$. 
Then the following statements are equivalent
\begin{description}
\item[(i)] $A$ has locally homogeneous singularities
\item[(ii)] There exist
a surjective
ring homomorphism $\rho:\; \C[[x_1, ... , x_n]] \arrow A$, 
where $\C[[x_1, ... , x_n]]$ is the ring of power series,   
and the ideal $\ker \rho$ is homogeneous in $\C[[x_1, ... , x_n]]$.
\end{description}

{\bf Proof:} Clear. \endproof

%%%%%%%%%%%%%%%%%%%%%%%%%%%%%%%%%%%%%%%%%%%%%%%%%%
\subsection{Hyperk\"ahler varieties with locally
homogeneous singularities}

\noindent
%%%%%%%%%%%%%%%%%%%%%%%%%%%%%%%%%%%%%%%%%%%%%%%%%%%%%%%%%%%%
\noindent\proposition\label{_comple_LHS<=>real_LHS_Proposition_}
Let $M$ be a complex variety, $M_\R$ the underlying real analytic
variety. Then $M_\R$ is a space with locally
homogeneous singularities (SLHS) if and only if $M$
is a space with locally
homogeneous singularities.

{\bf Proof:} This is \cite{_Verbitsky:Desingu_},
Proposition 4.6. \endproof

\hfill

%%%%%%%%%%%%%%%%%%%%%%%%%%%%%%%%%%%%%%%%%%%%%%%%%%%%%%%%%%%%%%%%%%%%%%
\corollary \label{_hype_SLHS_for_diff_indu_c_str_Corollary_}
\cite{_Verbitsky:Desingu_}
Let $M$ be a hyperk\"ahler variety, $I_1$, $I_2$ induced complex
structures. Then $(M, I_1)$ is a space with locally
homogeneous singularities if and only is $(M, I_2)$ is 
SLHS.

{\bf Proof:} The real analytic variety underlying 
 $(M, I_1)$ coinsides with that underlying 
 $(M, I_2)$. Applying \ref{_comple_LHS<=>real_LHS_Proposition_},
we immediately 
obtain \ref{_hype_SLHS_for_diff_indu_c_str_Corollary_}.
\endproof

\hfill

%%%%%%%%%%%%%%%%%%%%%%%%%%%%%%%%%%%%%%%%
\definition
Let $M$ be a hyperk\"ahler or hypercomplex variety. Then $M$ is called 
a space with locally homogeneous singularities (SLHS) if 
the underlying real analytic variety is SLHS
or, equivalently, for some induced complex structure
$I$ the $(M, I)$ is SLHS.

\hfill

Some of the canonical examples of hyperk\"ahler varieties
are spaces with locally homogeneous singularities {\it per se}. 
For instance, it is easy to prove the following theorem:

\hfill

%%%%%%%%%%%%%%%%%%%%%%%%%%%%%%%%%%%%%%%%
\theorem 
(\cite{_Verbitsky:Desingu_}, Theorem 4.9)
Let $M$ be a compact hyperk\"ahler manifold, $I$ an induced
complex structure and $B$ a stable holomorphic bundle over $(M, I)$.
Assume that $c_1(B)$, $c_2(B)$ are $SU(2)$-invariant, with respect
to the standard action of $SU(2)$ on $H^*(M)$. 
Let $\Def(B)$ be a reduction
of a deformation space of stable holomorphic 
structures on $B$, which is a hyperk\"ahler variety by 
\ref{_hyperho_defo_hyperka_Theorem_}. Then 
$\Def(B)$ is a space with locally homogeneous singularities (SLHS).

\endproof

\hfill

However, for the other examples of hyperk\"ahler varieties, 
there is no easy {\it ad hoc} way to show that they are SLHS. The main
aim of this paper, however, is to prove that every hypercomplex variety
is SLHS (see Subsection \ref{_main_resu_Subsection_}).

%%%%%%%%%%%%%%%%%%%%%%%%%%%%%%%%%%%%%%%%%%%%%%%%%%%%%%%%%%%%
\subsection{Desingularization of hypercomplex varieties which are SLHS}
\label{_desingu_for_SLHS_Subsection_}
%%%%%%%%%%%%%%%%%%%%%%%%%%%%%%%%%%%%%%%%%%%%%%%%%%%%%%%%%%%%

For hypercomplex varieties which are SLHS, we have a complete list of
possible singularities (\cite{_Verbitsky:Desingu_}; see also
\ref{_singula_stru_Theorem_}). 
This makes it possible to desingularize
every hypercomplex (or hyperk\"ahler) variety in a natural way. 
The present paper shows that {\em every} hypercomplex variety
is SLHS, thus extending the results of \cite{_Verbitsky:Desingu_}
to all hypercomplex varieties. For the benefit of the
reader, we relate in this Subsection 
the main results of \cite{_Verbitsky:Desingu_}. We don't use 
these results in the rest of the article, so the reader is free
to skip this Subsection.

\hfill

Here is the theorem describing the shape of singularities.

\hfill

%%%%%%%%%%%%%%%%%%%%%%%%%%%%%%%%%%%%%%%%%%%%%%%%%%
\theorem\label{_singula_stru_Theorem_}
Let $M$ be a hypercomplex variety, and $I$ an induced 
complex structure. Assume that $M$ is SLHS. Then, for each 
point $x\in M$, there exists a neighbourhood  $U$ of $x\in (M,I)$,
which is isomorphic to $B \cap\left(\bigcup_i L_i\right)$, where $B$ is an open
ball in $\C^n$ and $\bigcup_i L_i$  is a union of planes 
$L_i \in \C^n$ passing through $0\in \C^n$.
In particular, the normalization of $(M,I)$ is smooth.

{\bf Proof:} See Corollary 5.3 of \cite{_Verbitsky:Desingu_}.
\endproof

\hfill

Here is the desingularization theorem.

\hfill

%%%%%%%%%%%%%%%%%%%%%%%%%%%%%%%%%%%%%%%%
\theorem \label{_desingu_Theorem_}
(\cite{_Verbitsky:Desingu_}, Theorem 6.1)
Let $M$ be a hyperk\"ahler or a hypercomplex variety,
$I$ an induced complex structure.
Assume that $M$ is a space with 
locally homogeneous singularities.
Let \[ \widetilde{(M, I)}\stackrel n\arrow (M,I)\] 
be the normalization of
$(M,I)$. Then $\widetilde{(M, I)}$ is smooth and
has a natural hyperk\"ahler (respectively, hypercomplex) 
structure $\c H$, such that the associated
map $n:\; \widetilde{(M, I)} \arrow (M,I)$ agrees with $\c H$.
Moreover, the hyperk\"ahler (hypercomplex)
manifold $\tilde M:= \widetilde{(M, I)}$
is independent from the choice of induced complex structure $I$.

\endproof

%%%%%%%%%%%%%%%%%%%%%%%%%%%%%%%%%%%%%%%%%%%%%%%%%%
\subsection{The main result: every hypercomplex variety is 
SLHS}
\label{_main_resu_Subsection_}
%%%%%%%%%%%%%%%%%%%%%%%%%%%%%%%%%%%%%%%%%%%%%%%%%%

The proof of the following theorem is given in
Sections
\ref{_comple_with_au_Section_}--\ref{_homogeni_on_hype_Section_}.

%%%%%%%%%%%%%%%%%%%%%%%%%%%%%%%%%%%%%%%%
\theorem\label{_hyperco_SLHS_Theorem_}
(the main result of this paper)
Let $M$ be a hypercomplex variety. Then
$M$ is a space with locally homogeneous singularities
(SLHS).

\hfill

\ref{_hyperco_SLHS_Theorem_} has the following immediate 
corollary.

\hfill

%%%%%%%%%%%%%%%%%%%%%%%%%%%%%%%%%%%%%%%%%%%%%%%%%%
\corollary \label{_hyperco_desingu_Corollary_}
Let $M$ be a hypercomplex or a hyperk\"ahler variety. Then 
\ref{_singula_stru_Theorem_} (a theorem describing the shape
of the singularities of $M$) and \ref{_desingu_Theorem_}
(desingularization theorem) hold.

\endproof

%%%%%%%%%%%%%%%%%%%%%%%%%%%%%%%%%%%%%%%%

\section{Complete rings with automorphisms}
\label{_comple_with_au_Section_}

%%%%%%%%%%%%%%%%%%%%%%%%%%%%%%%%%%%%%%%%

%%%%%%%%%%%%%%%%%%%%%%%%%%%%%%%%%%%%%%%%
\definition \label{_homogeni_automo_Definition_}
Let $A$ be a local Noetherian ring over $\C$, equipped with an
automorphism
$e:\; A \arrow A$. Let $\goth m$ be a maximal ideal of $A$.
Assume that $e$ acts on $\goth m /\goth m^2$ as a multiplication 
by $\lambda\in \C$, $|\lambda|< 1$. Then $e$ is called {\bf a
homogenizing automorphism of $A$}.

\hfill

The aim of the present section is to prove the following statement.

%%%%%%%%%%%%%%%%%%%%%%%%%%%%%%%%%%%%%%%%%%%%%%%%%%%%%%
\proposition \label{_homogeni_LHS_Proposition_}
Let $A$ be a complete Noetherian ring over $\C$, equipped with a
homogenizing authomorphism $e:\; A \arrow A$. Then there exist
a surjective 
ring homomorphism $\rho:\; \C[[x_1, ... , x_n]] \arrow A$, 
such that the ideal
$\ker \rho$ is homogeneous in $\C[[x_1, ... , x_n]]$.
In particular, $A$ has locally homogeneous 
singularities (\ref{_locally_homo_coord_Claim_}).

\hfill

This statement is well known. A reader who knows its proof
should skip the rest of this section.

\hfill

%%%%%%%%%%%%%%%%%%%%%%%%%%%%%%%%%%%%%%%%%%%%%%%%%%%%%%%%%%%%
\proposition \label{_homogeni_auto_then_basis_Proposition_}
Let $A$ be a complete Noetherian ring over $\C$, equipped with a
homogenizing authomorphism $e:\; A \arrow A$. Then there exist
a system of ring elements
\[ 
    f_1 , ..., f_m \in \goth m, \ \ m = \dim_\C\goth m /\goth m^2,
\]
which generate $\goth m /\goth m^2$, and such that $e(f_i) = \lambda f_i$.

\hfill

{\bf Proof:}
Let $\underline a\in\goth m /\goth m^2$.
Let $a\in \goth m$ be a representative of $\underline a$ in $\goth m$.
To prove \ref{_homogeni_auto_then_basis_Proposition_}
it suffices to find $c \in \goth m^2$, such that 
$e(a-c) = \lambda a -\lambda c$. Thus, we need to solve an equation
\begin{equation}\label{_a_through_a_Equation_}
 \lambda c - e(c) = e(a) - \lambda(a). 
\end{equation}
Let $r:= e(a)-\lambda a$. Clearly, $r\in \goth m ^2$.
A solution of \eqref{_a_through_a_Equation_}
is provided by the following lemma.

\hfill

%%%%%%%%%%%%%%%%%%%%%%%%%%%%%%%%%%%%%%%%%%%%%%%%%%
\lemma \label{_e-lambda_invertible_Lemma_}
In assumptions of \ref{_homogeni_auto_then_basis_Proposition_},
let $r\in \goth m^2$. Then, the equation
\begin{equation}\label{_finding_eigen_Equation_}
e(c) - \lambda c = r
\end{equation}
has a unique solution $c \in \goth m^2$.

\hfill

{\bf Proof:} We need to show that the operator
$P:= (e-\lambda)\restrict{\goth m^2}$
is invertible. Consider the adic filtration 
$\goth m^2 \subset \goth m^3 \subset ...$ on $\goth m^2$.
Clearly, $P$ preserves this filtration. Since $\goth m^2$
is complete with respect to the adic filtration,
it suffices to show that $P$ is invertible on the
successive quotients. The quotient $\goth m^2/\goth m^i$ is 
finite-dimensional, so to show that $P$ is invertible it suffices
to calculate the eigenvalues. Since $e$ is an automorphism,
restriction of $e$ to $\goth m^i/\goth m^{i-1}$ is a multiplication
by $\lambda^i$. Thus, the eigenvalues of $e$ on $\goth m^2/\goth m^i$ 
range from $\lambda^2$ to $\lambda^{i-1}$. Since $|\lambda|>|\lambda|^2$,
all eigenvalues of $P\restrict{\goth m^2/\goth m^i}$ are 
non-zero and the restriction of $P$ to $\goth m^2/\goth m^i$ is invertible.
This proves \ref{_e-lambda_invertible_Lemma_}.
$\blacksquare$

\hfill

{\bf The proof of \ref{_homogeni_LHS_Proposition_}.}
Consider the map \[ \rho:\; \C[[x_1, ... x_m]] \arrow A,\ \ 
\rho(x_i) = f_i,\] where $f_1, ... , f_m$ is the system of functions
constructed in \ref{_homogeni_auto_then_basis_Proposition_}.
By Nakayama, $\rho$
is surjective. 

Let $e_\lambda:\; \C[[x_1, ... x_m]] \arrow\C[[x_1, ... x_m]] $ be the 
automorphism
mapping $x_i $ to $\lambda x_i$. Then, the diagram
\[\begin{CD} \C[[x_1, ... x_m]] @>{\rho}>> A \\
		@V{e_\lambda}VV @VV{e}V\\
		\C[[x_1, ... x_m]] @>{\rho}>> A 
\end{CD}
\]
is by construction commutative.
Therefore, the ideal $I= \ker \rho$ is preserved by $e_\lambda$.
Clearly, every $e_\lambda$-preserved 
ideal $I\subset \C[[x_1, ... x_m]]$ is homogeneous.
\ref{_homogeni_LHS_Proposition_} is proven. \endproof

%%%%%%%%%%%%%%%%%%%%%%%%%%%%%%%%%%%%%%%%%%%%%%%%%%%%%%%%%%%%%%%%%%%%%%

\section[Authomorphisms of local rings of holomorphic functions
on hyperk\"ahler varieties]
{Authomorphisms of local rings of holomorphic functions
on hyperk\"ahler varieties}
\label{_homogeni_on_hype_Section_}

%%%%%%%%%%%%%%%%%%%%%%%%%%%%%%%%%%%%%%%%%%%%%%%%%%%%%%%%%%%%%%%%%%%%%%

Let $M$ be a hypercomplex variety, $x\in M$ a point,
$I$ an induced complex structure. Let $A_I:= \hat \calo_x(M,I)$
be the adic completion of the local ring $\calo_x(M,I)$ of $x$-germs
of holomorphic functions on the complex variety $(M,I)$.
Clearly, the sheaf ring of the antiholomorphic functions on $(M,I)$
coinsides with $\calo_x(M,-I)$. Thus, the corresponding completion
ring is $A_{-I}$. As in \cite{_Verbitsky:Desingu_},
Claim 2.1, we have the natural isomorphism of completions:
\begin{equation}\label{_co_ana_and_rea_isom_Equation_}
\widehat{A_I \otimes_\C A_{-I}} =\widehat{\calo_x(M_\R)\otimes_\R \C},
\end{equation}
where \[ \widehat{\calo_x(M_\R)\otimes_\R \C}\] is the $x$-completion
of the ring of germs of real analytic complex-valued functions on $M$.
Consider the natural quotient map \[ p:\;A_{-I}\arrow \C.\]
Denote the ring \[ \widehat{\calo_x(M_\R)\otimes_\R \C}\] by $A_\R$. Let
$i_I:\; A_I \hookrightarrow A_\R$ be the natural embedding 
\[ a \mapsto a\times 1\in\widehat{A_I \otimes_\C A_{-I}},\]
and $e_I:\; A_\R \arrow A_I$ be the natural epimorphism associated
with the surjective map 
\[ A_I \otimes_\C A_{-I} \arrow A_I,\ \  
   a\otimes b \mapsto a\otimes p(b),
\] 
where $a\in A_I$, $b\in A_{-I}$, and
\[ a\otimes b\in{A_I \otimes_\C A_{-I}}\subset A_\R.\]
For an induced complex structure $J$, we define
$A_J$, $A_{-J}$, $i_J$, $e_J$ likewise.

Let $\Psi_{I,J}:\; A_I \arrow A_I$ be the composition

\[ A_I \stackrel {i_I}\arrow A_\R\stackrel {e_J}\arrow A_J
   \stackrel {i_J}\arrow A_\R\stackrel {e_I}\arrow A_I.
\]
Clearly, for $I=J$, the ring morphism $\Psi_{I,J}$ is identity,
and for $I=-J$, $\Psi_{I,J}$ is an augmentation map.

\hfill

%%%%%%%%%%%%%%%%%%%%%%%%%%%%%%%%%%%%%%%%%%%%%%%%%%%%%%%%%%%%
\proposition \label{_homogenizing_Proposition_}
Let $M$ be a singular hyperk\"ahler variety, $x\in M$ a point,
$I$, $J$ be induced complex structures, such that $I\neq J$ and
$I\neq -J$. Consider the map  $\Psi_{I,J}:\; A_I \arrow A_I$
defined as above. Then $\Psi_{I,J}$ is a homogenizing
automorphism of $A_I$.\footnote{For the definition
of a homogenizing automorphism, see \ref{_homogeni_automo_Definition_}.}

\hfill

{\bf Proof:}
Let $d\Psi$ be differential of $\Psi_{I,J}$, that is, 
the restriction of $\Psi_{I,J}$ to $\goth m/\goth m^2$,
where $\goth m$ is the maximal ideal of $A_I$. 
By Nakayama, to prove that
$\Psi_{I,J}$ is an automorphism it suffices to show that $d\Psi$
is invertible. To prove  that $\Psi_{I,J}$ is homogenizing, 
we have to show that $d\Psi$ is a multiplication by a complex
number $\lambda$, $|\lambda|<1$. As usually, we denote the real analytic
variety underlying $M$ by $M_\R$. 
Let $T_I$, $T_J$, $\underline {T}_\R$ be the Zariski 
tangent spaces to $(M,I)$, $(M,J)$ and $M_\R$, respectively, 
in $x\in M$. Consider the complexification
  $T_\R:= \underline {T}_\R\otimes \C$, which is a 
Zariski tangent space to the local ring $A_\R$.
To compute $d\Psi:\; T_I \arrow T_I$, we need
to compute the differentials of $e_I$, $e_J$, $i_I$, $i_J$,
i. e., the restrictions of the homomorphisms
$e_I$, $e_J$, $i_I$, $i_J$ to the Zariski tangent spaces
$T_I$, $T_J$, $T_\R$.
Denote these differentials by $de_I$, $de_J$, $di_I$, $di_J$.

\hfill

%%%%%%%%%%%%%%%%%%%%%%%%%%%%%%%%%%%%%%%%%%%%%%%%%%
\lemma \label{_i_e_through_Hodge_Lemma_}
Let $M$ be a hyperk\"ahler variety, $M_\R$ the associated real analytic
variety, $x\in M$ a point. Consider the space $T_\R := T_x (M_\R)\otimes
\C$. For an induced complex structure $I$, consider the Hodge decomposition
$T_\R= T^{1,0}_I \oplus  T^{0,1}_I$. In our previous notation,
$T_I^{1,0}$ is $T_I$. Then, $di_I$ is the natural embedding
of $T_I = T_I^{1,0}$ to $T_\R$, and $de_I$ is the natural
projection of $T_\R= T^{1,0}_I \oplus  T^{0,1}_I$ to
$T_I^{1,0}=T_I$.

{\bf Proof:} Clear. \endproof

\hfill

We are able now to describe the map $d\Psi:\; T_I \arrow T_I$
in terms of the quaternion action. Recall that the space $T_I$
is equipped with a natural $\R$-linear
quaternionic action. For each quaternionic
linear space $\underline V$ and each quaternion $I$, $I^2=-1$, $I$ defines a
complex structure in $\underline V$. Such a complex structure is called
{\bf induced by the quaternionic structure}. 

\hfill

%%%%%%%%%%%%%%%%%%%%%%%%%%%%%%%%%%%%%%%%%%%%%%%%%%
\lemma \label{_Psi_through_quate_Lemma_}
Let $\underline V$ be a space with quaternion action, and 
$V:= \underline V \otimes \C$ its complexification.
For each induced complex structure $I\in {\Bbb H}$,
consider the Hodge decomposition $V:= V_I^{1,0} \oplus V_I^{0,1}$.
For an induced complex structures $I, J\in \Bbb H$,
let $\Phi_{I,J}(V)$ be a composition of the natural embeddings
and projections
\[ 
   V_I^{1,0} \arrow V \arrow  V_J^{1,0} \arrow V \arrow  V_I^{1,0}.
\]
Using the natural identification $\underline V \cong  V_I^{1,0}$,
we consider $\Phi_{I,J}(V)$ as an $\R$-linear
automorphism of the space $\underline V$.
Then, applying the operator $\Phi_{I,J}(V)$  to the 
quaternionic space $T_I$, we obtain the operator $d\Psi$
defined above.

{\bf Proof:} Follows from \ref{_i_e_through_Hodge_Lemma_}
\endproof

\hfill

As we have seen, to prove \ref{_homogenizing_Proposition_}
it suffices to show that $d\Psi$ is a multiplication by a non-zero
complex number $\lambda$, $|\lambda| < 1$.
Thus, the proof of \ref{_homogenizing_Proposition_} is finished with the
following lemma.

\hfill

%%%%%%%%%%%%%%%%%%%%%%%%%%%%%%%%%%%%%%%%%%%%%%%%%%
\lemma\label{_compu_of_Psi_for_qua_Lemma_}
In assumptions of \ref{_Psi_through_quate_Lemma_}, consider the 
map \[ \Phi_{I,J}(V):\; V_I^{1,0} \arrow V_I^{1,0}.\] Then $\Phi_{I,J}(V)$
is a multiplication by a complex number $\lambda$. 
Moreover, $\lambda$ is a non-zero number unless $I=-J$,
and $|\lambda|< 1$ unless $I=J$.

\hfill

{\bf Proof:}
Let $\underline V= \oplus \underline V_i$ be a decomposition of $V$ into
a direct sum of $\Bbb H$-linear spaces. Then, the operator $\Phi_{I,J}(V)$
can also be decomposed: $\Phi_{I,J}(V) = \oplus \Phi_{I,J}(V_i)$.
Thus, to prove \ref{_compu_of_Psi_for_qua_Lemma_} it suffices
to assume that $\dim_{\Bbb H} \underline V=1$.
Therefore, we may identify $\underline V$ with the space
$\Bbb H$, equipped with the right action of quaternion
algebra on itself.

Consider the left action of $\Bbb H$ on $\underline V = \Bbb H$.
This action commutes with the right action of $\Bbb H$ on $\underline V$.
Consider the corresponding action 
\[ 
   \rho:\; SU(2) \arrow \End(\underline V)
\] of the group of unitary
quaternions ${\Bbb H}^{un}=SU(2)\subset \Bbb H$ on $\underline V$. 
Since $\rho$ commutes with the quaternion action, 
$\rho$ preserves $V^{1,0}_I \subset V$, for every
induced complex structure $I$. By the same token, for each $g\in SU(2)$,
the endomorphism $\rho(g)\in \End(V^{1,0}_I)$
commutes with $\Phi_{I,J}(V)$.

Consider the 2-dimensional $\C$-vector space $V^{1,0}_I$ 
as a representation of $SU(2)$. Clearly, $V^{1,0}_I$ 
is an irreducible representation. Thus, by Schur's lemma,
the automorphism $\Phi_{I,J}(V)\in \End(V^{1,0}_I))$ 
is a multiplication by a complex constant $\lambda$.
The estimation $0< |\lambda| < 1$ follows from the 
following elementary argument. The composition
$i_I \circ e_J$ applied to a vector 
$v\in V_I^{1,0}$ is a projection of $v$ to $V_J^{1,0}$
along $V_J^{0,1}$. Consider the natural Euclidean metric
on $V = \Bbb H$. Clearly, the decomposition 
$V = V_J^{1,0}\oplus V_J^{0,1}$ is orthogonal.
Thus, the composition $i_I \circ e_J$ is an orthogonal
projection of $v\in V_I^{1,0}$ to $V_J^{1,0}$.
Similarly, the composition $i_J \circ e_I$ is an orthogonal
projection of $v\in V_J^{1,0}$ to $V_I^{1,0}$.
Thus, the map $\Phi_{I,J}(V)$ is an orthogonal projection 
from $V_I^{1,0}$ to $V_J^{1,0}$ and back to $V_I^{1,0}$.
Such a composition always decreases a length of vectors,
unless $V_I^{1,0}$ coincides with $V_J^{1,0}$, in which
case $I=J$. Also, unless $V_I^{1,0} =  V_J^{0,1}$,
$\Phi_{I,J}(V)$ is non-zero; in the later case, $I = -J$.
\ref{_homogenizing_Proposition_}
is proven. This finishes the proof of 
\ref{_hyperco_SLHS_Theorem_}. \endproof

\hfill

{\bf Acknowledgements:} 
 A nice version of the proof
of \ref{_e-lambda_invertible_Lemma_} was suggested by Roma 
Bez\-ru\-kav\-ni\-kov.
I am grateful to A. Beilinson, R. Bez\-ru\-kav\-ni\-kov,
P. Deligne, D. Kaledin, D. Kazhdan, M. Kontsevich, V. Lunts,
T. Pantev and S.-T. Yau for enlightening discussions.

\end{document}